\documentclass[12pt]{amsart}
\usepackage{amsthm, amssymb}
\usepackage{amsfonts, amscd}
\usepackage{epsfig,multicol}
\theoremstyle{plain} \numberwithin{equation}{section}
\newtheorem{thm}{Theorem}[section]

\newtheorem{lem}{Lemma}[section]

\theoremstyle{definition}
\newtheorem{defn}[thm]{Definition}

\DeclareMathOperator{\Aut}{Aut}

\begin{document}
\title[Small Covers over Prisms]
{\large \bf Small Covers over Prisms}
\author[Mingzhong Cai, Xin Chen and Zhi L\"u]{Mingzhong Cai, Xin Chen and Zhi L\"u}
\footnote[0]{
\endgraf
 {\bf 2000AMS Classification:} 57S10, 57M60, 14M25, 52B70, 05C10.
 \endgraf
 Supported by grants from NSFC (No. 10371020 and No. 10671034).}
\address{School of Mathematical Sciences, Fudan University, Shanghai,
200433, P.R. China. Current address: Department of Mathematics,
Cornell University, Ithaca, NY, 14852, USA.}
\email{mc456@cornell.edu}
\address{Department of Statistics, University of California, Riverside,
CA, 92521, USA. } \email{xchen004@student.ucr.edu.}

\address{Institute of Mathematics, School of Mathematical Sciences, Fudan University, Shanghai,
200433, P.R. China.} \email{zlu@fudan.edu.cn}

\date{}
\begin{abstract} In this paper we calculate the number of equivariant diffeomorphism classes
of small covers over a prism.
\end{abstract}
\maketitle

\section{Introduction}

The notion of {\em small covers} is first introduced by Davis and
 Januszkiewicz \cite{dj}, where a small cover is a smooth closed
 manifold $M^n$ with a locally standard $({\Bbb Z}_2)^n$-action
 such that its orbit space is a simple convex polytope. This gives a direct connection between
 equivariant topology and combinatorics. As shown
 in \cite{dj}, all small covers over a simple convex polytope
 $P^n$ correspond to all characteristic functions (in this paper
 we call them $({\Bbb Z}_2)^n$-colorings) defined on all facets
 (codimension-one faces) of $P^n$. However, two small covers over
 $P^n$ which correspond to different $({\Bbb Z}_2)^n$-colorings
 may be equivariantly diffeomorphic. In \cite{lm}, L\"u and Masuda
 showed that the equivariant diffeomorphism class of a small cover over $P^n$
  agrees  with the equivalence class of its corresponding $({\Bbb
 Z}_2)^n$-coloring under the action of automorphism group of
 $P^n$. Therefore, the number of equivariant diffeomorphism
 classes of small covers over a certain polytope can be
 interpreted as that of equivalence classes of $({\Bbb
 Z}_2)^n$-colorings, and in particular, the latter one is enumerable.
 When $n=2$, a recursive formula for this number is obtained in
 \cite{lm}.  However, there is no instant formula to determine this
 number if $n>2$. In \cite{gs}, Garrison and Scott used a computer
 program to enumerate the number of  homeomorphism classes of all small covers
 over a dodecahedron, and the program only needs a few changes if
 we want to know the number of its equivariant diffeomorphism
 classes.

 \vskip .2cm

 Generally, the equivalence classes of $({\Bbb
 Z}_2)^n$-colorings of a polytope $P^n$ can be considered as the orbits
 of all $({\Bbb
 Z}_2)^n$-colorings under the action of automorphism group of
 $P^n$. Thus, Burnside Lemma (or Cauchy-Frobenius Lemma) can be applied to determine the
 number of the orbits
 of all $({\Bbb
 Z}_2)^n$-colorings under the action of automorphism group of
 $P^n$, which is just the average value of numbers $|\Lambda_g|$,
 where $|\Lambda_g|$ equals to the number of $({\Bbb
 Z}_2)^n$-colorings which are fixed by an automorphism $g$ of
 $P^n$. This leads us to give a formula for the number
 of equivariant diffeomorphism  classes of all small covers over a
 prism (see Theorem~\ref{n}), and in particular, the Euler's totient function is also involved.

\vskip .2cm

The arrangement of this paper is as follows. In Section 2 we
review  the basic theory  about small covers and Burnside Lemma,
and calculate the automorphism group of face poset of a prism. In
Section 3 we determine the number of all colorings on a prism, so
that in Section 4 we can obtain a calculation formula of the
number of equivariant diffeomorphism classes of all small covers
over a prism, and by using a computer program, the first 8 numbers
are given.

\section{Preliminaries}

\subsection{Small covers and colorings}

\vskip .2cm An $n$-dimensional convex polytope $P^n$ is said to be
{\em simple} if exactly $n$ faces of codimension one meet at each
of its vertices.   A closed $n$-manifold $M^n$ is said to be a
{\em small cover} if it admits an effective $({\Bbb
Z}_2)^n$-action, which is locally isomorphic to the standard
action of $({\Bbb Z}_2)^n$ on ${\Bbb R}^n$, and the orbit space of
the action is a simple convex polytope $P^n$.

\vskip .2cm Suppose that $\pi:M^n\longrightarrow P^n$ is a small
cover over a simple convex polytope $P^n$. Let
$\mathcal{F}(P^n)=\{F_1,...,F_\ell\}$ be the set of
codimension-one faces (facets) of $P^n$. Then there are $\ell$
connected submanifolds $M_1,...,M_\ell$ determined by $\pi$ and
$F_i$ (i.e., $M_i=\pi^{-1}(F_i)$). Each  submanifold $M_i$  is
fixed pointwise by the ${\Bbb Z}_2$-subgroup $G_i$ of  $({\Bbb
Z}_2)^n$, so that each facet $F_i$ corresponds to the ${\Bbb
Z}_2$-subgroup $G_i$. Obviously,  such the ${\Bbb Z}_2$-subgroup
$G_i$ actually agrees with an element $\upsilon_i$ in $({\Bbb
Z}_2)^n$ as a vector space.
 For each face $F$ of
codimension $u$, since $P^n$ is simple, there are $u$ facets
$F_{i_1},...,F_{i_u}$  such that $F=F_{i_1}\cap\cdots \cap
F_{i_u}$. Then, the corresponding characteristic submanifolds
$M_{i_1},...,M_{i_u}$ intersect transversally in the
$(n-u)$-dimensional submanifold $\pi^{-1}(F)$, and the isotropy
subgroup $G_F$ of $\pi^{-1}(F)$ is a  subtorus of rank $u$ and is
generated by $G_{i_1},...,G_{i_u}$ (or is determined by
$\upsilon_{i_1},...,\upsilon_{i_u}$ in $({\Bbb Z}_2)^n$).
 Thus,  this actually gives a characteristic
function (see \cite{dj})
$$\lambda:\mathcal{F}(P^n)\longrightarrow ({\Bbb Z}_2)^n$$
defined by $\lambda(F_i)=\upsilon_i$ such that for any face
$F=F_{i_1}\cap\cdots \cap F_{i_u}$ of $P^n$,
$\lambda(F_{i_1}),...,\lambda(F_{i_u})$ are linearly independent
in $({\Bbb Z}_2)^n$. If  we regard each nonzero vector of $({\Bbb
Z}_2)^n$ as being a {\em color}, then the characteristic function
$\lambda$ means  that each facet is colored by a color. Thus, we
also call $\lambda$ a {\em $({\Bbb Z}_2)^n$-coloring} on $P^n$
here.
 \vskip .2cm
Davis and Januszkiewicz \cite{dj} gave a reconstruction process of
$M^n$ by using the $({\Bbb Z}_2)^n$-coloring $\lambda$ and the
product bundle $({\Bbb Z}_2)^n\times P^n$ over $P^n$, so that all
small covers over $P^n$ are classified in terms of all $({\Bbb
Z}_2)^n$-colorings on ${\mathcal{F}}(P^n)$.   By $\Lambda(P^n)$ we
denote  the set of all $({\Bbb Z}_2)^n$-colorings on $P^n$. Then
we have

\begin{thm}
[Davis-Januszkiewicz] Let $\pi:M^n\longrightarrow P^n$ be a small
cover over a simple convex polytope $P^n$. Then all small covers
over $P^n$ are given by $\{M(\lambda)\vert \lambda\in
\Lambda(P^n)\}$.
\end{thm}

  \noindent {\bf Remark.} Generally speaking, we cannot make sure
that there always exist $({\Bbb Z}_2)^n$-colorings over a simple
convex polytope $P^n$ when $n\geq 4$.
    For example, see [DJ, Nonexamples 1.22]. However,
    the Four Color Theorem makes sure that every 3-dimensional simple
 convex polytope  always admits $({\Bbb Z}_2)^3$-colorings.

 \vskip .2cm

There is a natural action of $\text{GL}(n,{\Bbb Z}_2)$ on
$\Lambda(P^n)$ defined by the correspondence
$\lambda\longmapsto\sigma\circ\lambda$, and it is easy to see that
such an action on $\Lambda(P^n)$ is free. Without loss of
generality, we assume that  $F_1,...,F_n$ of $\mathcal{F}(P^n)$
meet at one vertex $p$ of $P^n$. Let $e_1,...,e_n$ be the standard
basis of $({\Bbb Z}_2)^n$. Write $A(P^n)=\{\lambda\in
\Lambda(P^n)\big| \lambda(F_i)=e_i, i=1,...,n\}$. Then we have

\begin{lem}\label{f}
$|\Lambda(P^n)|=|A(P^n)|\times |GL(n,{\Bbb Z}_2)|.$
\end{lem}

Note that we know from \cite{ab} that $|GL(n,{\Bbb
Z}_2)|=\prod_{k=1}^n(2^n-2^{k-1})$.

 \subsection{Automorphisms of polytopes and classification of
 small covers} Each point of a simple convex polytope $P^n$ has
a neighborhood which is affine isomorphic to an open subset of the
positive cone ${\Bbb R}_{\geq 0}^n$, so $P^n$ is an
$n$-dimensional manifold with corners (see \cite{d}).  An
automorphism of $P^n$ is a self-diffeomorphism of $P^n$ as a
manifold with corners, and by $\Aut(P^n)$ we denote the group of
automorphisms of $P^n$. All faces of $P^n$ forms a {\em poset}
(i.e., a partially ordered set by inclusion). An automorphism of
${\mathcal{F}}(P^n)$ is a bijection from ${\mathcal{F}}(P^n)$ to
itself which preserves the poset structure of all faces of $P^n$,
and by $\Aut({\mathcal{F}}(P^n))$ we denote the group of
automorphisms of ${\mathcal{F}}(P^n)$. Each automorphism of
$\Aut(P^n)$ naturally induces  an automorphism of
${\mathcal{F}}(P^n)$. It is well-known (see \cite{bp} or \cite{z})
that two simple convex polytopes are combinatorially equivalent
(i.e., their corresponding face posets are isomorphic) if and only
if they are diffeomorphic as manifolds with corners. Thus, the
natural homomorphism
$\phi:\Aut(P^n)\longrightarrow\Aut({\mathcal{F}}(P^n))$ is
surjective.

\begin{defn} Two $({\Bbb Z}_2)^n$-colorings $\lambda_1$ and
$\lambda_2$ in $\Lambda(P^n)$ are said to be equivalent if there
exists an automorphism $h\in\Aut({\mathcal{F}}(P^n))$ such that
$\lambda_1=\lambda_2\circ h$.
\end{defn}

Every small cover induces a $({\Bbb Z}_2)^n$-coloring on its orbit
polytope. However, it is possible that two $({\Bbb
Z}_2)^n$-colorings on a simple convex polytope  rebuild two
equivariantly diffeomorphic small covers. The following theorem
shows that the  equivalence of $({\Bbb Z}_2)^n$-colorings exactly
determines small covers up to equivariant diffeomorphism.

\begin{thm} \label{eq}
Two small covers over a simple convex polytope $P^n$ are
equivariantly diffemorphic if and only if their corresponding
$({\Bbb Z}_2)^n$-colorings on $P^n$ are equivalent.
\end{thm}

Readers may find the proof of Theorem~\ref{eq} in \cite{lm}, which
we omit here. Thus we may conclude that the equivalence classes of
$({\Bbb Z}_2)^n$-colorings on $P^n$ bijectively correspond to the
equivariant diffeomorphism classes of small covers over $P^n$.

\subsection{Burnside Lemma and automorphism group of a prism}
The equivalence classes of $({\Bbb Z}_2)^n$-colorings on $P^n$ can
be naturally considered as orbits of $\Lambda(P^n)$ under the
action of $\Aut({\mathcal{F}}(P^n))$. The famous Burnside Lemma
(or Cauchy-Frobenius Lemma) is essential in the enumeration of the
number of orbits.

\vskip .2cm

\noindent {\bf Burnside Lemma.} {\em Let $G$ be a finite group
acting on a set $X$. Then the number of orbits of $X$ under the
action of $G$ equals to ${1\over{|G|}}\sum_{g\in G}|X_g|$, where
$X_g=\{x\in X| gx=x\}$. }

\vskip .2cm

Burnside Lemma suggests that, in order to determine the number of
the equivalence classes of $({\Bbb Z}_2)^n$-colorings on $P^n$, we
need to understand the structure of $\Aut({\mathcal{F}}(P^n))$. As
stated in Section 1, we shall particularly be concerned with the
case in which $P^n$ is a prism.

\vskip .2cm

Let $P^3(m)$ be a $m$-sided prism (i.e., the product of a
$m$-polygon $P^2(m)$ and the interval $[0,1]$, or two $m$-polygons
joined by a belt of $m$ squares).

\begin{lem} \label{aut}
When $m\not=4$, the automorphism group
$\Aut({\mathcal{F}}(P^3(m)))$ is isomorphic to
$\mathcal{D}_{2m}\times {\Bbb Z}_2$, where $\mathcal{D}_{2m}$ is
the dihedral group of order $2m$. When $m=4$,
$\Aut({\mathcal{F}}(P^3(m)))$ is isomorphic to the direct product
$\mathcal{S}_4\times {\Bbb Z}_2$, where $\mathcal{S}_4$ is the
symmetric group of order 4.
\end{lem}
\begin{proof} All sided facets of $P^3(m)$ are 4-polygons, and the
top and bottom facets are two $m$-polygons.  If $m\not=4$,
obviously there are automorphisms of $P^3(m)$ under which  all
sided facets admit an action of dihedral group of order $2m$, and
there is also an automorphism of $P^3(m)$ such that  the top and
bottom facets is interchanged so they both admit an action of
${\Bbb Z}_2$. Since any one of  all sided facets cannot be mapped
to the top facet or bottom facet under the automorphisms of
$P^3(m)$, we have that the automorphism group
$\Aut({\mathcal{F}}(P^3(m)))$ is just isomorphic to the direct
product $\mathcal{D}_{2m}\times {\Bbb Z}_2$.

\vskip .2cm If $m=4$, then all facets of $P^3(4)$ consists of six
4-polygons so $P^3(4)$ is a 3-cube. A 3-cube has the same
automorphism group as an octahedron since a 3-cube is just the
dual of an octahedron. Obviously, the automorphism group of an
octahedron contains a symmetric group $\mathcal{S}_4$ of order 4
since there is exactly one such automorphism for each permutation
of the four pairs of opposite sides of the octahedron. In
addition, it is easy to see that an octahedron also admits a
reflection automorphism which is different from any one of
$\mathcal{S}_4$. Actually, such an automorphism reflects  two
opposite vertices of the octahedron. Thus,
$\Aut({\mathcal{F}}(P^3(4)))$ is isomorphic to the direct product
$\mathcal{S}_4\times {\Bbb Z}_2$.
\end{proof}

\noindent {\bf Remark.} It is not difficult to see that all
automorphisms of $P^3(4)$ contain those automorphisms that combine
a reflection and a rotation. Thus, $\Aut({\mathcal{F}}(P^3(4)))$
has three versions of $\mathcal{D}_8\times {\Bbb Z}_2$ as
subgroups.

\vskip .2cm

 By $s_1$ and $s_2$ we denote the top and bottom facets
of $P^3(m)$ respectively, and by $a_1,..., a_m$ we denote all
sided facets (i.e., 4-polygons) of $P^3(m)$ in their general
order. Let $x, y, z$ be three automorphisms of
$\Aut({\mathcal{F}}(P^3(m)))$ with the following properties
respectively:

$(1)\ \ x(a_i)=a_{i+1} (i=1,2,..., m-1), x(a_m)=a_1, x(s_j)=s_j,
j=1, 2;$

$(2)\ \ y(a_i)=a_{m+1-i} (i=1,2,..., m), y(s_j)=s_j, j=1, 2;$

$(3)\ \ z(a_i)=a_{i} (i=1,2,..., m),  z(s_1)=s_2, z(s_2)=s_1.$

\noindent Then, when $m\not=4$, all automorphisms of
$\Aut({\mathcal{F}}(P^3(m)))$ can be written in a simple form as
follows: \begin{equation} \label {gr} x^uy^vz^w, \ u\in{\Bbb Z}_m,
\ v,w \in{\Bbb Z}_2 \end{equation} with $x^m=y^2=z^2=1$,
$x^uy=yx^{m-u}$, and $x^uy^vz=zx^uy^v$.

\section{Colorings on a prism}

This section is devoted to calculating the number of $({\Bbb
Z}_2)^3$-colorings on a $m$-sided prism $P^3(m)$.

 \begin{thm} \label{t1}
Let $a, b, c$ be the functions from ${\Bbb N}$ to ${\Bbb N}$ with
the following properties:

$1)\ \  a(j)=2a(j-1)+8a(j-2)$ with $a(1)=1, a(2)=2$;

$2)\ \  b(j)=b(j-1)+4b(j-2)$ with $b(1)=b(2)=1$;

$3)\ \  c(j)=2c(j-1)+4c(j-2)-6c(j-3)-3c(j-4)+4c(j-5)$ with
$c(1)=c(2)=1, c(3)=3, c(4)=7, c(5)=17$.  Then the number of
$({\Bbb Z}_2)^3$-colorings on a $m$-sided prism $P^3(m)$ is
$$|\Lambda(P^3(m)|=168[a(m-1)+2b(m-1)+c(m-1)].$$
 \end{thm}

 \begin{proof}
$({\Bbb Z}_2)^3$ contains seven nonzero elements (or seven colors)
$e_1, e_2, e_3, e_1+e_2, e_1+e_3, e_2+e_3, e_1+e_2+e_3$ where
$e_1, e_2, e_3$ form a standard basis of $({\Bbb Z}_2)^3$. Set
$$A(m)=\{\lambda\in \Lambda(P^3(m))\big|
\lambda(s_1)=e_1, \lambda(a_1)=e_2,\lambda(a_2)=e_3\}.$$ Then, by
Lemma~\ref{f}, we have that $|\Lambda(P^3(m))|=|A(m)|\times
|GL(3,{\Bbb Z}_2)|=168|A(m)|$. Write
$$\begin{cases}
A_1(m)=\{\lambda\in A(m)\big|
\lambda(s_2)=e_1\} \\
A_2(m)=\{\lambda\in A(m)\big|
\lambda(s_2)=e_1+e_2\} \\
A_3(m)=\{\lambda\in A(m)\big|
\lambda(s_2)=e_1+e_3\} \\
A_4(m)=\{\lambda\in A(m)\big| \lambda(s_2)=e_1+e_2+e_3\}.
\end{cases}$$
By the definition of $({\Bbb Z}_2)^3$-colorings, it is easy to see
that $|A(m)|=\sum_{i=1}^4|A_i(m)|$. Then, our argument proceeds as
follows.

\vskip .2cm

{\bf Case 1.} Calculation of $|A_1(m)|$.

\vskip .2cm Take a coloring $\lambda$ in $A_1(m)$, by the
definition of $({\Bbb Z}_2)^3$-colorings, we easily see that
$\lambda(a_m)$ has four possible values $e_3, e_1+e_3, e_2+e_3,
e_1+e_2+e_3$, and it is also so even if $\lambda(a_1)=e_1+e_2$
since $\lambda(s_1)=\lambda(s_2)=e_1$. Set $A_1^0(m)=\{\lambda\in
A_1(m)| \lambda(a_{m-1})=e_2 \text{ or } e_1+e_2\}$ and
$A_1^1(m)=A_1(m)\setminus A_1^0(m)$. Take a coloring $\lambda$ in
$A_1^0(m)$. Then $\lambda(a_{m-1})=e_2$ or $e_1+e_2$, and so the
possible values of $\lambda(a_{m-2})$ are $e_3, e_1+e_3, e_2+e_3,
e_1+e_2+e_3$. In this case, we see that the values of $\lambda$
restricted to facets $a_{m-1}$ and $a_{m}$ have only eight
possible choices. Thus, $|A_1^0(m)|$ is just eight times of
$|A_1(m-2)|$. Take a coloring $\lambda$ in $A_1^1(m)$. Then the
possible values of $\lambda(a_{m-1})$ are $e_3, e_1+e_3, e_2+e_3,
e_1+e_2+e_3$.  In this case, if we fix any one of four possible
values of $\lambda(a_{m-1})$, then it is easy to see that
$\lambda(a_{m})$ has only two possible values. Thus, $|A_1^1(m)|$
is just two times of $|A_1(m-1)|$. Further, we have that
$$|A_1(m)|=2|A_1(m-1)|+8|A_1(m-2)|.$$
A direct observation  shows that when $m=2$, $|A_1(m)|=1$, and
when $m=3$, $|A_1(m)|=2$. Thus, we have that $|A_1(m)|=a(m-1)$.

\vskip .2cm

{\bf Case 2.} Calculation of $|A_2(m)|$.

\vskip .2cm Similarly to the case 1, set $A_2^0(m)=\{\lambda\in
A_2(m)| \lambda(a_{m-1})=e_2 \}$ and $A_2^1(m)=A_2(m)\setminus
A_2^0(m)$. Take a coloring $\lambda$ in $A_2^0(m)$, we have that
each of both $\lambda(a_{m-2})$ and $\lambda(a_m)$ has four
possible values $e_3, e_1+e_3, e_2+e_3, e_1+e_2+e_3$, so
$|A_2^0(m)|=4|A_2(m-2)|$; take a coloring $\lambda$ in $A_2^1(m)$,
we then have that $\lambda(a_{m-1})$ has four possible values
$e_3, e_1+e_3, e_2+e_3, e_1+e_2+e_3$ but $\lambda(a_{m})$ has only
one possible value whichever of four possible values of
$\lambda(a_{m-1})$ is chosen, so $|A_2^1(m)|=|A_2(m-1)|$. Also, we
easily see that $|A_2(2)|=|A_2(3)|=1$. Thus, $|A_2(m)|=b(m-1)$.

\vskip .2cm

{\bf Case 3.} Calculation of $|A_3(m)|$.

\vskip .2cm

If we interchange $e_2$ and $e_3$ in this case, then the problem
is reduced to the case 2, so $|A_3(m)|=b(m-1)$. \vskip .2cm

{\bf Case 4.} Calculation of $|A_4(m)|$.

\vskip .2cm

In this case, given a coloring $\lambda\in A_4(m)$, we have that
for any sided facet $a_i$,  $\lambda(a_i)$ cannot be $e_1$ and
$e_1+e_2+e_3$. Note that $\lambda(a_m)$ is equal to $e_3$ or
$e_2+e_3$ since $\lambda(s_2)=e_1+e_2+e_3$. Set
$A_4^0(m)=\{\lambda\in A_4(m)| \lambda(a_{m-1})=e_2\}$,
$A_4^1(m)=\{\lambda\in A_4(m)| \lambda(a_{m-1})=e_3 \text{ or }
e_2+e_3\}$, and $A_4^2(m)=\{\lambda\in A_4(m)|
\lambda(a_{m-1})=e_1+e_2 \text{ or } e_1+e_3\}$. Then
$|A_4(m)|=|A_4^0(m)|+|A_4^1(m)|+|A_4^2(m)|$. An easy argument
shows that $|A_4^0(m)|=2|A_4(m-2)|$ and $|A_4^1(m)|=|A_4(m-1)|$,
so
\begin{equation} \label{e1}
|A_4(m)|=|A_4(m-1)|+2|A_4(m-2)|+|A_4^2(m)|.
\end{equation} Set  $B(m)=\{\lambda\in A_4^2(m)|
\lambda(a_{m-2})=e_2+e_3\}$. Then it is easy to see that
\begin{equation}\label{e2}
|A_4^2(m)|=|A_4^2(m-1)|+|B(m)| \end{equation}  and
\begin{equation}\label{e3}
|B(m)|=2|A_4(m-4)|+2|A_4(m-5)|+|B(m-2)|+2|A_4^2(m-2)|.
\end{equation}
Combining equations (\ref{e1}), (\ref{e2}) and (\ref{e3}), we
obtain
$$|A_4(m)|=2|A_4(m-1)|+4|A_4(m-2)|-6|A_4(m-3)|-3|A_4(m-4)|+4|A_4(m-5)|.$$
A direct observation gives that $|A_4(2)|=|A_4(3)|=1, |A_4(4)|=3,
|A_4(5)|=7$, and $|A_4(6)|=17$. Thus, we have $|A_4(m)|=c(m-1)$.
 \end{proof}

 Note that an easy argument shows  $a(j)={1\over
6}[4^j-(-1)^j2^j]$. Also, using a different argument, we can
obtain another expression of $c(j)$, i.e.,
$c(j)=2c(j-1)+3c(j-2)-4c(j-3)+(-1)^{j}2$.

\section{The number of equivariant diffeomorphism classes}

Let $\varphi(n)$ be the value at $n\in{\Bbb N}$ of the Euler's
totient function, which is the number of positive integers that
are both less than $n$ and coprime to $n$. Specifically, if we
write $n=p_1^{l_1}\cdots p_r^{l_r}$ where the $p_i$ are distinct
primes, then
$\varphi(n)=(p_1^{l_1}-p_1^{l_1-1})\cdots(p_r^{l_r}-p_r^{l_r-1})$.
Note that $\varphi(1)=1$.  Let $\nu(m)$ denote the number of
$({\Bbb Z}_2)^3$-colorings on $P^3(m)$ with its top and bottom
facets having the same color, so $\nu(m)=168a(m-1)$. Let $E(m)$
denote the number of equivariant diffeomorphism classes of small
covers over $P^3(m)$, which is just the number of the equivalence
classes of all $({\Bbb Z}_2)^3$-colorings on $P^3(m)$.

\begin{thm} \label{n}
The number $E(m)$ is equal to
$$\begin{cases}
{1\over{4m}}\big\{\sum_{k>1, k|m}\varphi({m\over
k})[|\Lambda(P^3(k))|+\nu(k)]+21m\rho_1(m)+42m\rho_2(m)\big\} &
\text{ if
$m\not=4$}\\
259   & \text{ if $m=4$}
\end{cases}$$
where $\rho_1(m)$ is defined recursively as follows
$$\rho_1(m)=\begin{cases}
0 & \text{ $m$ is odd}\\
5 & \text{ $m=0$}\\
12 & \text{ $m=2$}\\
\rho_1(m-2)+4\rho_1(m-4) & \text{ otherwise}
\end{cases}$$
and $$\rho_2(m)=\begin{cases} 0 & \text{ $m$ is odd}\\
2^m & \text{ $m$ is even.}
\end{cases}$$
\end{thm}

\begin{proof}
 From Burnside Lemma and Lemma~\ref{aut}, we have that
$$E(m)=\begin{cases} {1\over
{4m}}\sum_{g\in \Aut({\mathcal{F}}(P^3(m)))}|\Lambda_g| & \text{
if $m\not=4$}\\
{1\over {48}}\sum_{g\in \Aut({\mathcal{F}}(P^3(4)))}|\Lambda_g| &
\text{ if $m=4$}
\end{cases}$$
where $\Lambda_g=\{\lambda\in \Lambda(P^3(m))|
\lambda=\lambda\circ g\}$.

\vskip .2cm  If $m\not=4$, then by (\ref{gr}) each automorphism
$g$ of $\Aut({\mathcal{F}}(P^3(m))$ can be written as $x^uy^vz^w$,
and the argument is divided into the following cases.

\vskip .2cm

{\bf Case 1:} $g=x^u$.

\vskip .2cm In this case, $g$ is just an automorphism which
rotates all sided facets of $P^3(m)$.  Let $k=\gcd(u, m)$. Then
all sided facets of $P^3(m)$ are divided into $k$ orbits under the
action of $g$, and each orbit contains ${m\over k}$ facets. Thus,
each $({\Bbb Z}_2)^3$-coloring of $\Lambda_g$ gives the same
coloring on all ${m\over k}$ facets of each orbit. This means that
if $k\not=1$, $|\Lambda_g|$  exactly equals to  the number of
$({\Bbb Z}_2)^3$-colorings on a $k$-sided prism, so
$|\Lambda_g|=|\Lambda(P^3(k))|$. If $k=1$, then all sided facets
have the same coloring, which is impossible by the definition of
$({\Bbb Z}_2)^3$-colorings. On the other hand, for every $k>1$,
there are exactly $\varphi({m\over k})$ automorphisms of the form
$x^u$, each of which  divides  all sided facets of $P^3(m)$ into
$k$ orbits. Thus, when $g$ is of the form $x^u$,
$$\sum_{g=x^u}|\Lambda_g|=\sum_{k>1, k|m}\varphi({m\over
k})|\Lambda(P^3(k))|.$$

\vskip .2cm

{\bf Case 2:} $g=x^uz$.

\vskip .2cm

Each such automorphism $g=x^uz$ gives an interchange between top
and bottom facets, so each $({\Bbb Z}_2)^3$-coloring of
$\Lambda_g$ gives the same coloring on top and bottom facets.
Combining the argument of the case 1, when $g$ is of the form
$x^uz$,
$$\sum_{g=x^uz}|\Lambda_g|=\sum_{k>1, k|m}\varphi({m\over
k})\nu(k).$$

\vskip .2cm

{\bf Case 3:} $g=x^uy$ or $x^uyz$ with $m$ odd.

\vskip .2cm

In this case, since $m$ is odd, each automorphism always reflects
at least two neighborly sided facets, so that the two neighborly
sided facets have the same coloring. But this contradicts the
definition of $({\Bbb Z}_2)^3$-colorings. Thus, for each such an
automorphism $g$,  $\Lambda_g$ is empty.

\vskip .2cm

{\bf Case 4:} $g=x^uy$ or $x^uyz$ with $u$ even and $m$ even.

\vskip .2cm

Let $l={{m-u-2}\over 2}$. Then it is easy to see that such an
automorphism in this case  gives an interchange between two sided
facets $a_l$ and $a_{l+1}$, so both facets $a_l$ and $a_{l+1}$
have the same coloring, exactly as in the case 3. Thus, in this
case $\Lambda_g$ is empty, too.

\vskip .2cm

{\bf Case 5:} $g=x^uy$ with $u$ odd and $m$ even.

\vskip .2cm

Since each automorphism $g=x^uy$ contains a reflection $y$ as its
factor and $u$ is odd, it becomes a reflection along a plane
through the center of some sided facet. Thus,  each coloring
$\lambda$ of $\Lambda_g$ is equivalent to coloring only ${m\over
2}+1$ neighborly sided facets and top and bottom facets of
$P^3(m)$. We shall show that for each $g=x^uy$, the number of all
colorings in $\Lambda_g$ is just
\begin{equation} \label{n1}
|\Lambda_g|=42\rho_1(m)+42\rho_2(m) \end{equation} where
$\rho_1(m)$ and $\rho_2(m)$ are stated as in Theorem~\ref{n}. It
is easy to see that there are exactly ${m\over 2}$ such
automorphisms $g=x^uy$ since $m$ is even and $u$ is odd, so
$$\sum_{g=x^uy}|\Lambda_g|=21m\rho_1(m)+21m\rho_2(m).$$
Now let us show the equality (\ref{n1}) as follows.
 \vskip .2cm

 Actually, the argument method of Case 1 of Theorem~\ref{t1} can
 still be carried out here. Also, it suffices to consider the case
 $g=x^{m-1}y$ (i.e., $g=yx)$ since other cases have no difference essentially.
 Set $X_1(m)=\{\lambda\in \Lambda_g|
 \lambda(s_1)\not=\lambda(s_2)\}$ and $X_2(m)=\Lambda_g\setminus
 X_1(m)$. Then, by $X_1^0(m)$ we denote the set $\{\lambda\in
 X_1(m)| \lambda(a_m)\not=\lambda(s_1)+\lambda(a_2),
 \lambda(a_m)\not=\lambda(s_2)+\lambda(a_2), \lambda(a_m)\not=\lambda(a_2)\}$, and by $X_1^1(m)$
 denote $X_1(m)\setminus X_1^0(m)$. Similarly to the argument of Case 1 of
 Theorem~\ref{t1}, we have that $|X_1^0(m)|=|X_1(m-2)|$ and
 $|X_1^1(m)|=4|X_1(m-4)|$ with initial values $X_1(0)=5\times6\times 7 $ and
 $X_1(2)=12\times 6\times 7$. Thus, $|X_1(m)|=42\rho_1(m)$. For
 $X_2(m)$, in a similar way we may obtain $|X_2(m)|=4|X_2(m-2)|$ with $X_2(0)=42$,
 which is exactly $42\rho_2(m)$.

  \vskip .2cm

{\bf Case 6:} $g=x^uyz$ with $u$ odd and $m$ even.

\vskip .2cm

This case is the same as $X_2(m)$ of Case 5. Thus,
$\sum_{g=x^uyz}|\Lambda_g|=21m\rho_2(m)$.

\vskip .2cm

Combining Cases 1-6, we complete the proof in the case $m\not=4$.

\vskip .2cm

If $m=4$, then $P^3(4)$ is a 3-cube. The automorphism group of a
3-cube has order 48, and contains $\mathcal{D}_8\times {\Bbb Z}_2$
of order 16 as a subgroup. As shown above, actually we have
determined the case of the action of $\mathcal{D}_8\times {\Bbb
Z}_2$ on $P^3(4)$. However, each of other 32 automorphisms of
$P^3(4)$ has no fixed coloring in $\Lambda(P^3(4))$ since it maps
top facet (or bottom facet) to a sided facet. Thus
$$E(4)={1\over{48}}\{\sum_{k=2,4}\varphi({4\over
k})[|\Lambda(P^3(k))|+\nu(k)]+84\rho_1(4)+168\rho_2(4)\}=259.$$
\end{proof}

\noindent {\bf Remark.} We may write a computer program to
calculate $E(m)$ by the formula in Theorem~\ref{n}. The first 8
numbers are stated as follows.
$$\begin{tabular}{|l|l|l|l|l|l|l|l|l|}
 \hline $m$ &
 3 & 4& 5& 6& 7&  8& 9&10\\
\hline $E(m)$ &  98&  259&  882&  4200&  9114& 35406&
 107086& 394632 \\
\hline
\end{tabular}$$

\end{document}